\documentclass{elsart3-1}

 \usepackage{graphicx}

\usepackage{amssymb,amsmath}
\usepackage{url}
\usepackage{subfigure}
\newenvironment{system}{\begin{equation} \left\{ \begin{array}{l}}{\end{array} \right.\end{equation}}

\usepackage[english,francais]{babel}

\newtheorem{e-proposition}[theorem]{Proposition}

\newtheorem{e-definition}[theorem]{Definition\rm}

\newtheorem{theoreme}{Th\'eor\`eme}[section]

\newtheorem{proposition}[theoreme]{Proposition}

\renewcommand{\theequation}{\arabic{equation}}
\def\og{\leavevmode\raise.3ex\hbox{$\scriptscriptstyle\langle\!\langle$~}}
\def\fg{\leavevmode\raise.3ex\hbox{~$\!\scriptscriptstyle\,\rangle\!\rangle$}}

\journal{the Acad\'emie des sciences}
\begin{document}
\centerline{}
\begin{frontmatter}


\selectlanguage{english}
\title{An analytical solution of Shallow Water system coupled to Exner equation}


\selectlanguage{english}
\author[authorlabel1]{C. Berthon},
\ead{christophe.berthon@math.univ-nantes.fr}
\author[authorlabel2]{S. Cordier},
\ead{stephane.cordier@math.cnrs.fr}
\author[authorlabel3]{O. Delestre}
\ead{Delestre@unice.fr} and
\author[authorlabel4]{M.H. Le}
\ead{mh.le@brgm.fr}
\address[authorlabel1]{Laboratoire de Mathematiques Jean Leray, Universit\'e de Nantes -- 2 rue de la Houssiniere -- 44322 Nantes}
\address[authorlabel2]{MAPMO UMR CNRS 6628, Universit\'e d'Orl\'eans, UFR Sciences, B\^atiment de
math\'ematiques -- 45067 Orl\'eans}
\address[authorlabel3]{Laboratoire de Math\'ematiques J.A. Dieudonn\'e \& Ecole Polytech Nice - Sophia, Universit\'e de Nice - Sophia Antipolis, Parc Valrose -- 06108 Nice}
\address[authorlabel4]{BRGM -- 3 avenue Claude Guillemin - B.P. 36009 -- 45060 Orl\'eans Orl\'eans}

\medskip
\begin{center}
{\small Received *****; accepted after revision +++++\\
Presented by }
\end{center}

\begin{abstract}
\selectlanguage{english}
In this paper, an exact smooth solution for the equations modeling the bedload transport of sediment in Shallow Water is presented. This solution is valid for a large family of sedimentation laws which are widely used in erosion modeling such as the Grass model or those of Meyer-Peter \& M\"uller. One of the main interest of this solution is the derivation of numerical benchmarks to valid the approximation methods. {\it To cite this article: A. Name1, A. Name2, C. R. Acad. Sci. Paris, Ser. I 340 (2005).}

\vskip 0.5\baselineskip

\selectlanguage{francais}
\noindent{\bf R\'esum\'e} \vskip 0.5\baselineskip \noindent
{\bf Une solution analytique du syst\`eme de Saint-Venant coupl\'e \`a l'\'equation d'Exner. }
Ce papier pr\'esente une solution analytique pour le syst\`eme mod\'elisant le transport de s\'ediments par le charriage. Cette solution est valable pour une grande famille de lois s\'edimentaires comme le mod\`ele de Grass ainsi que celui de Meyer-Peter \& M\"uller. Ce r\'esultat est utile pour la validation des sch\'emas num\'eriques.
{\it Pour citer cet article~: A. Name1, A. Name2, C. R. Acad. Sci.
Paris, Ser. I 340 (2005).}

\end{abstract}
\end{frontmatter}

\selectlanguage{francais}


\selectlanguage{english}
\section{Introduction}\label{intro}

Soil erosion is a consequence of the movements of sediments due to mechanical actions of flows. In the context of bedload transport, a mass conservation law, also called Exner equation \cite{exner25}, is used to update the bed elevation. This equation is coupled with the shallow water equations describing the overland flows (see \cite{castro08sediment} and references therein) as follows:
\begin{align}
& \partial_t h + \partial_x (hu)=0, \label{eq01}\\
& \partial_t (hu)+\partial_x(hu^2+gh^2/2)+gh\partial_x z_b= 0, \label{eq02}\\
& \partial_t z_b + \partial_x q_b = 0, \label{eq03}
\end{align}
where $h$ is the water depth, $u$ the flow velocity, $z_b$ the thickness of sediment layer which can be modified by the fluid and $g$ the acceleration due to gravity. The variable $hu$ is also called water discharge and noted by $q$. Finally, $q_b$ is the volumetric bedload sediment transport rate. Its expressions are usually proposed for granular non-cohesive sediments and quantified empirically \cite{grass81sediment,meyer-peter48formulas,Cordier2011980}.

Many numerical schemes have been developed to solve system (\ref{eq01}-\ref{eq03}) (see \cite{Cordier2011980}
and references therein). The validation of such schemes by an analytical solution is a simple way to ensure their working. Nevertheless, analytical solutions are not proposed in the literature. Up to our knowledge, asymptotic solutions, derived by Hudson in \cite{Hudson2001}, are in general adopted to perform some comparisons with approximated solutions. The solutions are derived for Grass model \cite{grass81sediment}, {\it i.e} $q_b=A_g u^3$,  when the interaction constant $A_g$ is smaller than $10^{-2}$. In this paper, we propose a non obvious analytical solution in the steady state condition of flow.


\section{Solution of the equations}
In order to obtain an analytic solution, we consider $q_b$ as a function of the dimensionless bottom shear stress $\tau^*_b$ (see \cite{Cordier2011980}). By using the friction law of Darcy-Weisbach, $\tau^*_b$ is given by
\begin{equation*}\label{eq:tau}
\tau^*_b = \frac{fu^2}{8(s-1)gd_s},
\end{equation*}
where $f$ is the friction coefficient, $s=\rho_s/\rho$ the relative density of sediment in water and $d_s$ the diameter of sediment. The formul\ae~ of $q_b$ is usually expressed under the form
\begin{equation}\label{eq:qb}
q_b = \kappa(\tau^*_b-\tau^*_{cr})^p\sqrt{(s-1)gd^3_s},
\end{equation}
where $\tau^*_{cr}$ is the threshold for erosion, $\kappa$ an empirical coefficient and $p$ an exponent which is often fixed to $3/2$ in many applications. 
The expression \eqref{eq:qb} can be written in the simple form
 \begin{equation}\label{eq:qb simple}
 q_b = Au_e^{2p},
 \end{equation}
 where the effective velocity $u_e$ and the interaction coefficient $A$ are defined by
 \begin{system}\label{eq:notation}
 u_e^2 = u^2-u_{cr}^2,\\
 u_{cr}^2 = \tau^*_{cr}\Big[\dfrac{f}{8(s-1)gd}\Big]^{-1},\\
 A = \kappa\Big[\dfrac{f}{8(s-1)gd}\Big]^p\sqrt{(s-1)gd^3_s}.
 \end{system}

\noindent{\bf Remark.} The Grass model \cite{grass81sediment} is one of the simplest case by using $p= 3/2,~\tau^*_{cr} = 0$ and an empirical coefficient $A_g$ instead of $A$. The Meyer-Peter \& M\"uler model \cite{meyer-peter48formulas} is one of the most applied by using $p=3/2,~\kappa=8,~\tau^*_{cr} =0.047$. The following result is valid for all models rewriting in form (\ref{eq:qb simple}-\ref{eq:notation}).

\begin{proposition} \label{prop}
Assume that $q_b$ is defined by \eqref{eq:qb}. For a given uniform discharge $q$ such that $\tau^*_b>\tau^*_{cr}$, system (\ref{eq01}-\ref{eq03}) has the following analytical unsteady solution
\begin{system}\label{solution}
  u^2_e = \Big[\dfrac{\alpha x + \beta}{A}\Big]^{1/p},\\
  u = \sqrt{u_e^2+u_{cr}^2},~h=q/u,\\
  z_b^0 = -\dfrac{u^3+2gq}{2gu}+C,\\
  z_b = -\alpha t + z_b^0.
\end{system}
where $\alpha,~\beta,~C$ are constants and $A,~u_{cr}$ are defined by \eqref{eq:notation}.
\end{proposition}

\noindent{\bf Proof.} We are here concerned by the smooth solution. In view of the assumption $hu = q = \text{cst}$, equations (\ref{eq01}-\ref{eq03}) reduce to
\begin{align}
& \partial_t h = 0, \notag\\
& \partial_x(q^2/h)+gh\partial_x H =0, \label{eq05}\\
& \partial_t H +\partial_x q_b = 0, \label{eq06}
\end{align}
where $H=h+z$ is the free surface elevation. Differentiating equation \eqref{eq05} with respect to $t$ and then equation \eqref{eq06} with respect to $x$, we obtain
\begin{align}
& \partial_{xt} H = 0,\notag\\
& \partial^2_x q_b = 0\label{eq08}.
\end{align}
Note that we can write $q_b=q_b(h,q)$ to have $\partial_t q_b = \partial_h q_b\partial_t h+\partial_q q_b\partial_t q = 0$, so $q_b$ is not time-depending. Thank to \eqref{eq08}, the expression of $q_b$ is obtained under the form
\begin{equation}\label{eq09}
  q_b = \alpha x + \beta,
\end{equation}
where $\alpha$ and $\beta$ are constant. From \eqref{eq03}, we obtain $\partial_t z_b = -\partial_x q_b = -\alpha$ to write
\begin{equation}\label{eq:zb}
  z_b = -\alpha t + z_b^0(x).
\end{equation}
Moreover, from \eqref{eq:qb simple} we deduce the effective velocity as follows:
\begin{equation*}\label{eq:ue}
  u_e^2 = \Big[\dfrac{\alpha x + \beta}{A}\Big]^{1/p}.
\end{equation*}
Plugging \eqref{eq:zb} into the momentum equation \eqref{eq05} and using a direct calculation, we have
\begin{equation*}\label{eq:zb0}
\partial_x z_b^0 = \Big[\frac{q}{u^2}-\frac{u}{g}\Big]\partial_x u \Rightarrow z_b^0 = -\frac{u^3+2gq}{2gu}+C
\end{equation*}
which concludes the proof.

\noindent{\bf Remark.} As $h$ and $u$ are stationary, the initial condition of \eqref{solution} is $(h,u,z_b^0)$. Moreover, the solution $(h,u)$ applied to the Grass model is also an analytical solution of the Shallow Water Equations with the variable topography $z_b^0$. Concerning the shallow-water model, other solutions can be found in \cite{SWASHES}.

\section{Numerical experiments}
\begin{figure}
\begin{tabular}{cc}
  \includegraphics[angle=-90,width = 0.49 \textwidth]{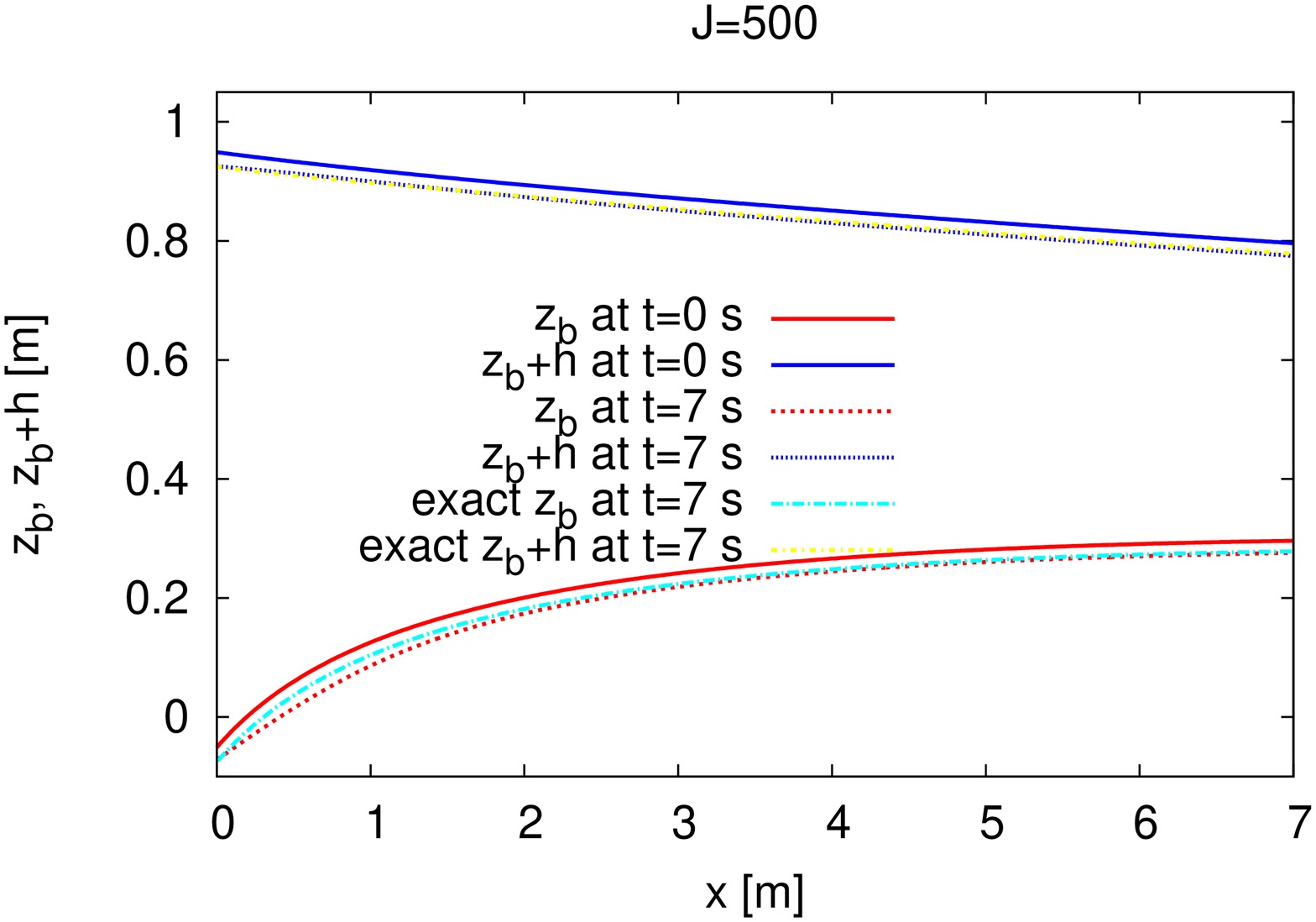} & \\[-5.9cm]
  &\includegraphics[angle=-90,width = 0.49 \textwidth]{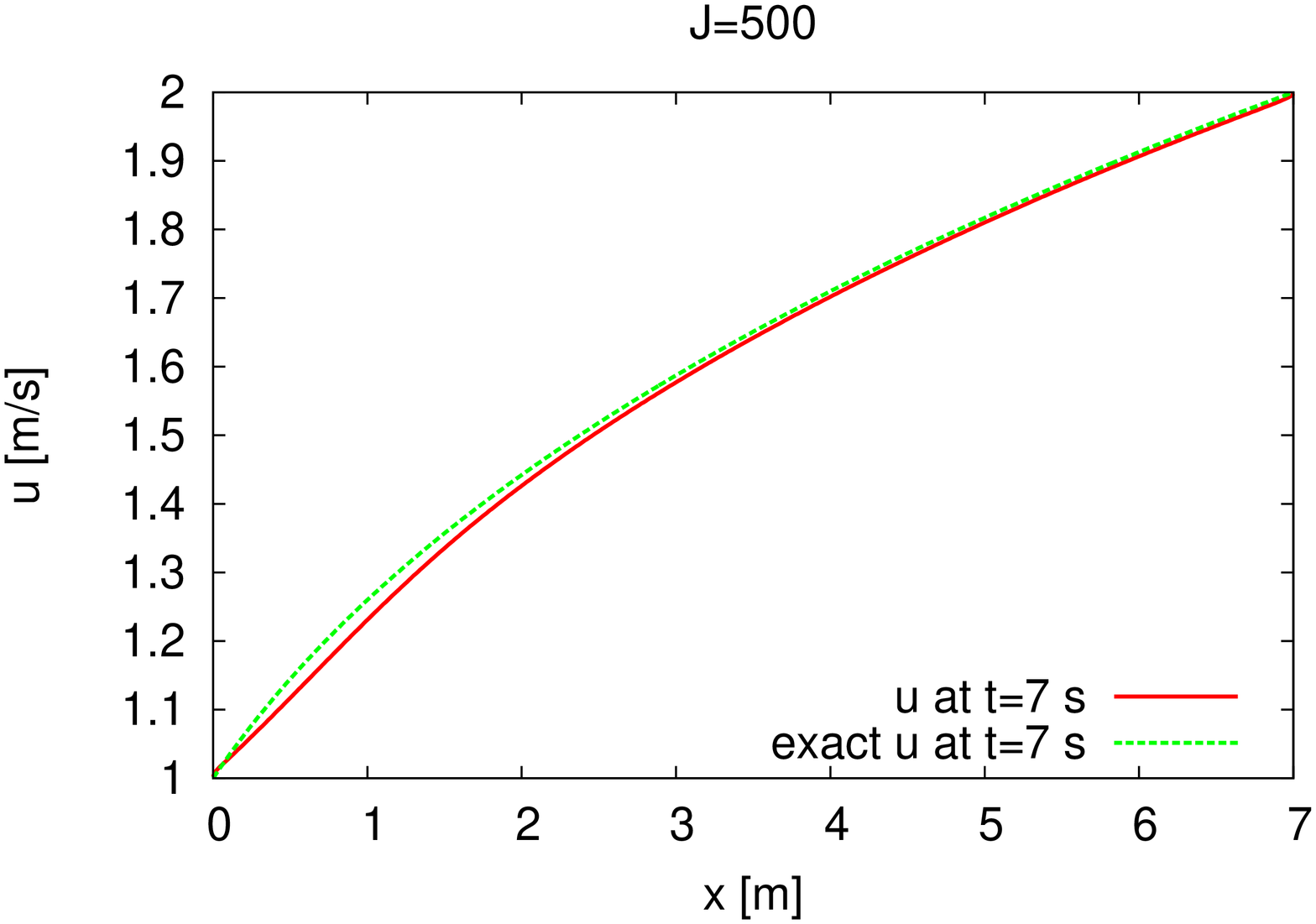}
\end{tabular}
\caption{Comparison between the exact solution and the relaxation method for : the water height and the topography (left) and
 the velocity (right).}
\label{fig-relaxation}
\end{figure}

In this section, we consider the analytical solution \eqref{solution} applied to the Grass model with $q=1$, $A_g=\alpha=\beta=0.005$ and $C=1$. A relaxation solver is applied to approximate the solution of the model. We will not give here the details of the relaxation solver (for details see \cite{Audusse2011}), but just the relaxation model for the equations (\ref{eq01}-\ref{eq03}). Thus, we solve the following relaxation system:
\begin{equation*}
\left\{\begin{array}{l}
 \partial_t h + \partial_x (hu)=0, \\
 \partial_t (hu)+\partial_x(hu^2+\pi)+gh\partial_x z_b= 0, \\
 \partial_t \pi + u\partial_x \pi + \frac{a^2}{h}\partial_x u = 0, \\
 \partial_t z_b + \partial_x q_r = 0, \\
 \partial_t q_r + \left(\frac{b^2}{h^2}-u^2\right)\partial_x z_b + 2u \partial_x q_r = 0,
\end{array}\right.
\end{equation*}
that is completed with $\pi=gh^2/2$ and $q_r=q_b$ at the equilibrium. 
 Figure \ref{fig-relaxation} presents the numerical result with $J=500$ space cells, a CFL fix condition of $1$ and $T = 7s$. We only notice little difference on the velocity, near the inflow boundary.

\begin{thebibliography}{00}



\bibitem{exner25} F. Exner. \"{U}ber die wechselwirkung zwischen wasser und geschiebe in fl\"{u}ssen,
Sitzungsber., Akad. Wissenschaften pt. IIa; 1925. Bd. 134.

\bibitem{castro08sediment} M.J.~Castro D\'{\i}az, E.D. Fern\'{a}ndez-Nieto, and A.M.Ferreiro.
 Sediment transport models in shallow water equations and numerical approach by high order finite volume methods.
 \emph{Computers \& Fluids}, 37\penalty0 (3):\penalty0 299--316, March 2008.

\bibitem{grass81sediment} A.J. Grass. Sediment transport by waves and currents.
 \emph{SERC London Cent. Mar. Technol}, Report No. FL29, 1981.

\bibitem{meyer-peter48formulas} E.~Meyer-Peter and R.~M{\"u}ller. Formulas for bed-load transport. In
 \emph{2nd meeting IAHSR, Stockholm, Sweden}, pages 1--26, 1948.


\bibitem{Cordier2011980} S.~Cordier, M.H. Le, and T.~Morales de~Luna.
 Bedload transport in shallow water models: Why splitting (may) fail, how hyperbolicity (can) help.
 \emph{Advances in Water Resources}, 34\penalty0 (8):\penalty0 980 -- 989, 2011.

\bibitem{Hudson2001} J.~Hudson. \emph{Numerical technics for morphodynamic modelling}.
 PhD thesis, University of Whiteknights, 2001.

\bibitem{SWASHES} O.~Delestre, C.~Lucas, P.-A.~Ksinant, F.~Darboux, C.~Laguerre, T.N.T.~Vo, F.~James and S.~Cordier. SWASHES: a library of Shallow Water Analytic Solutions for Hydraulic and Environmental Studies (\emph{Submitted}), http://arxiv.org/abs/1110.0288

\bibitem{Audusse2011} E.~Audusse, C.~Chalons, O.~Delestre, N. Goutal, M.~Jodeau, J. Sainte-Marie and B.~Spinewine.
Sediment transport modelling : Three layer models and relaxation schemes.
\emph{In preparation : CEMRACS 2011}.



\end{thebibliography}
\end{document}